\documentclass[runningheads]{llncs}
\usepackage{graphicx}
\usepackage{multirow}
\usepackage{subcaption}
\usepackage{cite} 
\usepackage{stmaryrd}
\usepackage{amsmath}
\usepackage{amsfonts}

\SetSymbolFont{stmry}{bold}{U}{stmry}{m}{n}

\usepackage{ifthen}
\usepackage{algorithm}
\usepackage{algorithmic}

\newcommand{\mv}{\boldsymbol{m}}
\newcommand{\cov}{\boldsymbol{C}}

\newcommand{\x}{\boldsymbol{x}}
\newcommand{\y}{\boldsymbol{y}}
\newcommand{\z}{\boldsymbol{z}}
\newcommand{\s}{\boldsymbol{s}}

\newcommand{\E}{\mathbb{E}}

\newcommand{\R}{\mathbb{R}}
\newcommand{\X}{\mathcal{X}}

\renewcommand{\S}{\mathcal{S}}

\newcommand{\N}{\mathcal{N}}

\newcommand{\argmin}{\mathop{\rm arg~min}\limits}

\newcommand{\T}{\mathrm{T}}

\newcommand{\I}{\mathbf{I}}

\newcommand{\p}{\boldsymbol{p}}

\usepackage{xcolor}
\usepackage{ifthen}
\newcommand{\markupdraft}[2]{
    \ifthenelse{\equal{#1}{display}}{#2}{}
    \ifthenelse{\equal{#1}{color}}{\color{#2}}{}
}

\newcommand{\newcolored}[3][]{{\markupdraft{color}{#2}#3}
    \ifthenelse{\equal{#1}{}}{}{\markupdraft{display}{{\color{yellow!70!black}[#1]}}}}
\newcommand{\del}[2][]{{\markupdraft{display}{{\color{orange}[removed: ``#2''[#1]]}}}} 
\newcommand{\new}[2][]{\newcolored[#1]{blue}{#2}}
\newcommand{\nnew}[2][]{\newcolored[#1]{red}{#2}}

\renewcommand{\del}[2]{}  
\renewcommand{\markupdraft}[2]{}  

\newcommand{\rev}[1]{{#1}}
\newcommand{\revdel}[1]{}

\begin{document}
\title{CMA-ES for Discrete and Mixed-Variable Optimization on Sets of Points}
%
%

\author{
Kento Uchida\inst{1}
\and Ryoki Hamano\inst{2}
\and Masahiro Nomura\inst{2}
\and Shota Saito\inst{1,3}
\and Shinichi Shirakawa\inst{1}
}

\authorrunning{K. Uchida et al.}

\institute{Yokohama National University, Yokohama, Japan\\
\email{
uchida-kento-fz@ynu.ac.jp, 
saito-shota-bt@ynu.jp,
shirakawa-shinichi-bg@ynu.ac.jp}
\and 
CyberAgent, Inc., Shibuya, Japan\\
\email{
hamano\_ryoki\_xa@cyberagent.co.jp,
nomura\_masahiro@cyberagent.co.jp,
}
\and
SKILLUP NeXt Ltd., Chiyoda, Japan\\
}


%
\maketitle              
\begin{abstract}

Discrete and mixed-variable optimization problems have appeared in several real-world applications.
Most of the research on mixed-variable optimization considers a mixture of integer and continuous variables, and several integer handlings have been developed to inherit the optimization performance of the continuous optimization methods to mixed-integer optimization.
\new{In some applications, acceptable solutions are given by selecting possible points in the disjoint subspaces.}
This paper focuses on the optimization on sets of points and proposes an optimization method by extending the covariance matrix adaptation evolution strategy (CMA-ES), termed the CMA-ES on sets of points (CMA-ES-SoP).
The CMA-ES-SoP incorporates margin correction that maintains the generation probability of neighboring points to prevent premature convergence to a specific non-optimal point, which is an effective integer-handling technique for CMA-ES.
\new{In addition, because margin correction with a fixed margin value tends to increase the marginal probabilities for a portion of neighboring points more than necessary,}
the CMA-ES-SoP updates the target margin value adaptively to make the average of the marginal probabilities close to a predefined target probability.
Numerical simulations demonstrated that the CMA-ES-SoP successfully optimized the optimization problems on sets of points, whereas the naive CMA-ES failed to optimize them due to premature convergence.

\end{abstract}
%
%
%

\section{Introduction}


Mixed-variable optimization methods have been actively developed due to the significant demand in real-world applications.
Most of the existing works have focused on mixed-integer optimization problems that contain both continuous and integer variables.
One of the major approaches is applying integer handlings to powerful continuous optimization methods to address integer variables~\cite{cmaesim,cmaeswm,dxnesici}.
For example, the reference~\cite{cmaeswm} focused on the covariance matrix adaptation evolution strategy (CMA-ES)~\cite{hansen:1996:ec} and proposed the CMA-ES with margin by incorporating the margin correction.
\nnew{The margin correction uses a margin that is a lower bound of the marginal probabilities for integer variables and prevents premature convergence caused by the original update of the CMA-ES.}
These integer handlings consist of coordinate-wise operations for each integer variable.

\revdel{However, other types of mixed-variable optimization have often been appeared.}
\rev{However, because the integer handling assumes the set of possible values are given on grid space, they cannot be applied to other kinds of sets of possible values.}
For example, when optimizing the location for the construction of wind turbines~\cite{wind-turbine}, the user makes a \revdel{list}\rev{set} of possible locations \rev{(pairs of latitude and longitude)} and requests an optimizer to select the best location from the \revdel{list}\rev{set}.
In this case, the existing integer handling cannot be applied.
In addition, when optimizing both the location and forms of winds that are represented by continuous variables, this problem is formulated as a mixed-variable optimization problem.
We term this problem structure as an optimization problem on the sets of points, and we formulate this problem as an optimization problem on the search space consisting of several subspaces.
Each subspace contains multiple possible points where the objective function value can be computed.
In mixed-variable optimization, some of the subspace is treated as continuous space.
The optimization problems on the sets of points have been found in several real-world applications such as design optimization of vehicle~\cite{car-design,mazda-design} and facility layout optimization~\cite{facility:2007}.
\new{We note that the naive CMA-ES fails to optimize such optimization problems on the sets of points, which will be observed in our experimental results.}

\begin{figure}[t]
    \centering
    \begin{minipage}[b]{0.3\linewidth}
        \centering
        \includegraphics[scale=0.3]{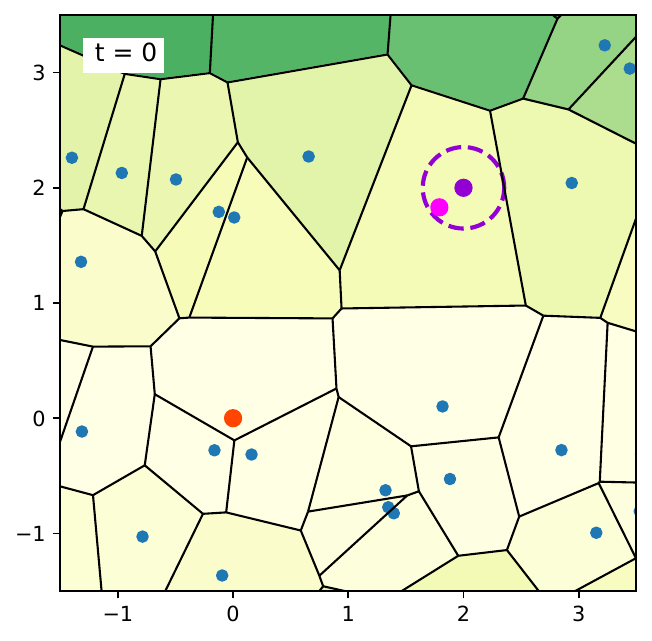}
    \end{minipage}
    \begin{minipage}[b]{0.3\linewidth}
        \centering
        \includegraphics[scale=0.3]{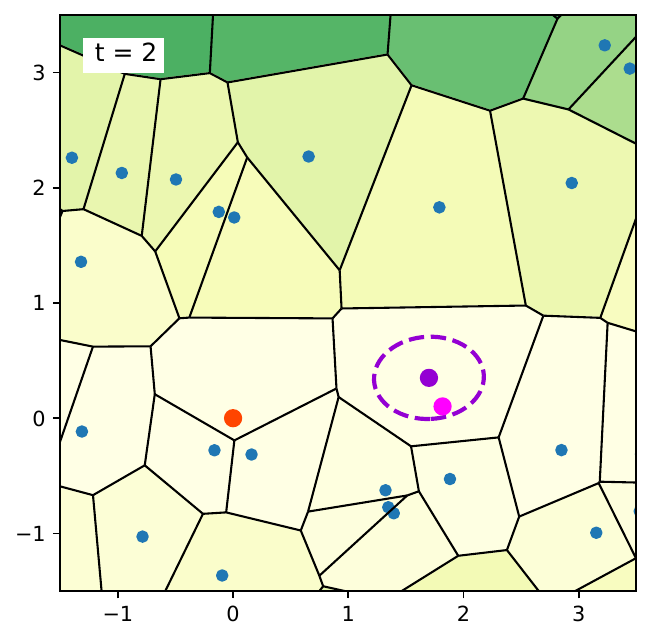}
    \end{minipage}
    \begin{minipage}[b]{0.3\linewidth}
        \centering
        \includegraphics[scale=0.3]{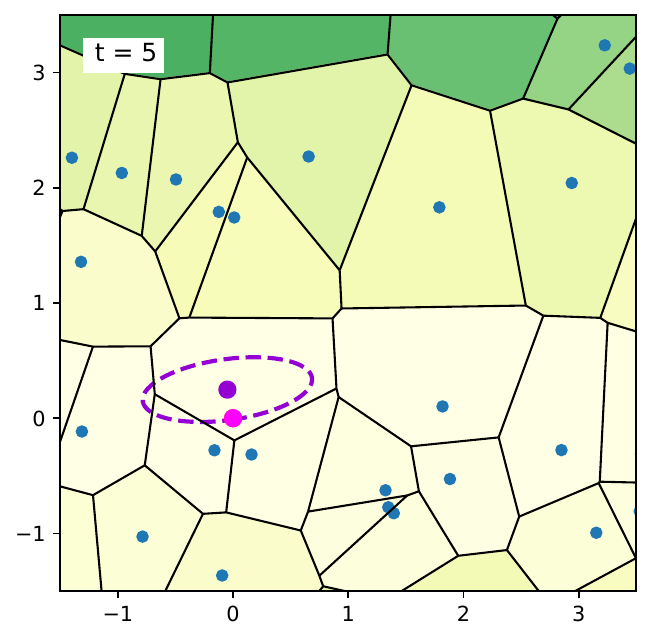}
    \end{minipage} 
    \caption{
        Illustration of the search space and the optimization process of CMA-ES-SoP on the two-dimensional Ellipsoid function. 
        The acceptable solutions are depicted as blue points. 
        The red and magenta points are the optimum and closest points, respectively.
        The CMA-ES-SoP uses the Voronoi diagram to encode the samples and adjust the margin value.
    }
    \label{fig:concept}
\end{figure}
%
%
%


In this paper, we tailor the CMA-ES for optimization on the sets of points and propose CMA-ES-SoP (CMA-ES on sets of points).
\new{Figure~\ref{fig:concept} shows the conceptual image of the optimization with the CMA-ES-SoP.}
The CMA-ES-SoP incorporates three handlings: sample encoding, margin correction, and margin adaptation.
In the encoding process, the samples generated from a multivariate Gaussian distribution are \revdel{mapped}\rev{projected} to the closest points \new{to the samples}.
In the margin correction, the covariance matrix is modified to maintain the marginal generation probability beyond \new{the mid-points between the mean vector and} neighboring points in the Voronoi diagram above the margin.
Finally, \new{to prevent an unnecessary increase of marginal probabilities after the margin adaptation}, the margin is adjusted so that the average of marginal probabilities is maintained close to the target value for the margin.

We evaluated the performance of CMA-ES-SoP using numerical simulations with benchmark functions.
In the experiment with discrete optimization on sets of points, the CMA-ES-SoP successfully optimized the benchmark functions with high probability, while the CMA-ES failed to optimize them.
In the experiment with mixed-variable optimization, the CMA-ES-SoP outperformed the CMA-ES in most functions, especially in high-dimensional problems.

\paragraph{Notations}
The functions $\Phi_{\mathrm{cdf}}: \R \to (0,1)$ and $\Phi_{\mathrm{ppf}}: (0,1) \to \R$ are the cumulative density function of the standard normal distribution $\mathcal{N}(0,1)$ and the inverse function of $\Phi_{\mathrm{cdf}}$ called the percentile point function, respectively.
We denote the concatenation of $n$ vectors $\boldsymbol{v}_1, \cdots, \boldsymbol{v}_n$ as $\textsc{Concat}(\boldsymbol{v}_1, \cdots, \boldsymbol{v}_n) = (\boldsymbol{v}_1^\T, \cdots, \boldsymbol{v}_n^\T)^\T$.

\section{Related Works}
\paragraph{Evolutionary algorithms for discrete optimization:}
Integer optimization, which is the optimization of integer variables, is a related problem \del{structure }{}to the optimization on sets of points.
Several evolutionary algorithms (EAs) have been developed for integer optimization.
The coordinate-wise mutation is a widely-used approach for integer optimization~\cite{Doerr:gecco:2016,Doerr:algorithmica:2018}.
Particularly, (1+1)-EA with self-adjusting mutation is a promising method with theoretical guarantee~\cite{Doerr:algorithmica:2018}.
Although there are several other approaches, including binary encoding~\cite{binaryenc:2002} and the probabilistic model-based approach~\cite{rEDA:2023}, there is no approach that can directly be applied to the optimization on sets of points.

The optimization problem with categorical variables is another topic for EAs~\cite{asng:2019,Berny:gecco:2021}.
As the points in the set can be treated as categories, these methods can be applied to the optimization on sets of points.
However, as the positional relationship between points is not addressed, \del{the optimization performance is improvable}{}\new{the optimization performance is limited}.
In addition, these EAs cannot deal with mixed-variable optimization problems containing both discrete and continuous variables.

\paragraph{Integer handling for mixed-integer optimization:}
Another related work is optimization methods for mixed-integer optimization.
Several powerful optimization methods have been developed by introducing integer handlings to powerful continuous optimization methods such as the CMA-ES.
The study~\cite{cmaesim} injects the integer mutation vector into the generation process of candidate solutions in the CMA-ES.
The CMA-ES with margin~\cite{cmaeswm} incorporates the margin correction to maintain the marginal probabilities associated with the integer variables \del{above the margin}{}\new{above a certain value}.
The DX-NES-ICI~\cite{dxnesici} \del{proposes the leap of }{}\new{leaps} the elements of the mean vector corresponding to integer variables to overcome the performance deterioration of CMA-ES with margin when the evaluation value is more influenced by continuous variables than integer variables.
Although these are powerful methods for mixed-integer optimization problems, they do not adequately handle scenarios involving a mix of continuous variables and variables on a set of points.

\section{Baseline Algorithm: CMA-ES}
\label{sec:cmaes}
CMA-ES~\cite{hansen:1996:ec} is a powerful black-box optimization method on continuous space.
The CMA-ES employs a multivariate Gaussian distribution as a sampling distribution of the candidate solutions and updates the distribution parameters to generate better solutions.
The multivariate Gaussian distribution is parameterized by the mean vector $\mv \in \R^N$, the step-size $\sigma \in \R_{>0}$, and the covariance matrix $\cov \in \R^{N \times N}$. The CMA-ES also employs two evolution paths $\p_c \in \R^N$ and $\p_\sigma \in \R^N$ that are initialized to zero vectors.

The update procedure of the CMA-ES is as follows.
First, the CMA-ES generates $\lambda$ solutions $\x^{\langle 1 \rangle}, \cdots, \x^{\langle \lambda \rangle}$ as
\begin{align}
    \z^{\langle i \rangle} \sim \N(\mathbf{0}, \I) 
    \enspace, \quad
    \y^{\langle i \rangle} = \sqrt{\cov^{(t)}} \z^{\langle i \rangle} 
    \enspace\text{, and}\quad
    \x^{\langle i \rangle} = \sigma^{(t)} \y^{\langle i \rangle} + \mv^{(t)} \enspace.
    \label{eq:cma:generation}
\end{align}
Subsequently, the generated solutions are evaluated on the objective function $f: \R^N \to \R$ to be optimized.
We denote the index of the $i$-th best solution as $i\!:\!\lambda$.

Next, the CMA-ES updates the evolution paths as 
\begin{align}
    \p_\sigma^{(t+1)} &= (1 - c_\sigma) \p_\sigma^{(t)} + \sqrt{ c_\sigma (2 - c_\sigma) \mu_\mathrm{eff} } \cdot \sum_{i = 1}^{\mu} w_i \z^{\langle i:\lambda \rangle} \\
    \p_c^{(t+1)} &= (1 - c_c) \p_c^{(t)} + h_\sigma^{(t+1)} \sqrt{ c_c (2 - c_c) \mu_\mathrm{eff} } \cdot \sum_{i = 1}^{\mu} w_i \y^{\langle i:\lambda \rangle}
    \enspace,
\end{align}
where $\{ w_i \}_{i=1}^\mu$ are predefined positive weights, $c_\sigma, c_c \in \R_{>0}$ are the accumulation rates of the evolution paths, and $\mu_\mathrm{eff} = (\sum_{i=1}^\mu w_i^2)^{-1}$ is the  variance effective selection mass.
The Heaviside function $h_\sigma^{(t+1)} \in \{0, 1\}$ becomes $h_\sigma^{(t+1)} = 1$ if and only if it satisfies:
\begin{align}
    \frac{\| \p_\sigma^{(t+1)} \| }{\sqrt{1 - (1 - c_\sigma)^{2 (t + 1)}}} < \left( 1.4 + \frac{2}{N+1} \right) \chi_N
    \enspace,
\end{align}
where $\chi_N = \sqrt{N} \left( 1 - \frac{1}{4N} + \frac{1}{21 N^2} \right)$ is the approximated value of the expectation $\E[ \| \N(\mathbf{0}, \I) \|]$. Otherwise, it becomes $h_\sigma^{(t+1)} = 0$.

Finally, the CMA-ES updates the distribution parameters of the multivariate Gaussian distribution as
\begin{align}
    \mv^{(t+1)} &= \mv^{(t)} + c_m \sum_{i = 1}^{\mu} w_i \left( \x^{\langle i:\lambda \rangle} - \mv^{(t)} \right)
    \label{eq:cma:update-mean} \\
    \sigma^{(t+1)} &= \sigma^{(t)} \exp \left( \frac{c_\sigma}{d_\sigma} \left( \frac{ \| \p_\sigma^{(t+1)} \| }{\chi_N} - 1 \right) \right)
    \label{eq:cma:update-stepsize}     
    \allowdisplaybreaks[3] \\
    \begin{split}
        \cov^{(t+1)} &= \left(1 + \delta (h_\sigma^{(t+1)}) \right) \cov^{(t)} + c_1 \left( \p_c^{(t+1)} \left( \p_c^{(t+1)} \right)^\T - \cov^{(t)} \right) \\
        &\hspace{70pt} + c_\mu \sum_{i = 1}^{\mu} w_i \left( \y^{\langle i:\lambda \rangle} \left( \y^{\langle i:\lambda \rangle} \right)^\T - \cov^{(t)} \right)
        \enspace,
    \end{split}
    \label{eq:cma:update-cov} 
\end{align}
where $c_m, c_1, c_\mu \in \R_{>0}$ are the learning rates, $d_\sigma \in \R_{>0}$ is the damping factor, and $\delta(h) = (1 - h) c_1 c_c (2 - c_c)$.
The CMA-ES has well-tuned default values for all hyperparameters. 
Details are available in the literature~\cite{hansen:2017:arxiv,hansen:2014:book}.

\section{Target Problem}
In this study, we consider the search space \del{$\X \in \R^N$}{}\new{$\X \subseteq \R^N$} consisting of $K$ subspaces as
\begin{align}
    \X = {\S}_1 \times \cdots \times {\S}_K \enspace.
\end{align}
\new{In discrete optimization problems on sets of points}, each subspaces ${\S}_k \subseteq \R^{N_k}$ is given by sets of $L_k$ points\new{, i.e., ${\S}_k = \{ {\s}_{k,1}, \cdots, {\s}_{k,L_k} \}$}\del{ $\{ {\s}_{k,1}, \cdots, {\s}_{k,L_k} \}$}{}, where $N = \sum_{k=1}^K N_k$.
\new{We assume that the sets of points are accessible for the optimization methods.}
The candidate solutions are constructed by selecting points from subspaces.
\nnew{In mixed-variable optimization problems, \del{on the other hand, }{}a part of search space is given by the continuous space, and other space is given by sets of points, i.e., $\X = {\S}_1 \times \cdots \times {\S}_{K-1} \times \R^{N_{\mathrm{co}}}$.}

Such problems can be found in applications of system design and manufacturing.
For example, in the design optimization of vehicle~\cite{car-design,mazda-design}, the optimization is sometimes performed by selecting available parts of machines. 
In the optimization of position for constructing wind turbines~\cite{wind-turbine}, the potential places can be listed in advance, and \del{then we can optimize the position by selecting}{}\new{the optimal place is selected} from the listed possible places. 
Additionally, when optimizing both the position and form of the wing, which is controlled by continuous variables, the problem becomes a mixed-variables optimization problem.
Similar problem structures can also be found in facility layout problems~\cite{facility:2007} and location-routing problems~\cite{lrp:2014}.

\new{
As well on the mixed-integer optimization problems, the naive CMA-ES usually fails to optimize the optimization on the sets of points.
This is because the original update rule of the CMA-ES leads to premature convergence around non-optimal points, as observed in our experimental results shown in Section~\ref{sec:exp}.
}

\section{Proposed Method: CMA-ES-SoP}
We propose a novel variant of the CMA-ES for optimization on sets of points, termed CMA-ES-SoP.
The CMA-ES-SoP encodes the samples generated from the multivariate Gaussian distribution to obtain the candidate solution on the search space.
After updating the distribution parameters, the CMA-ES-SoP corrects the covariance matrix $\cov^{(t+1)}$ to maintain the generation probability \del{of}{}\new{beyond the mid-points from the} neighboring points above the margin value $\alpha^{(t)}$.
\new{In addition, we adapt the margin value $\alpha_k^{(t)} \in \R_{>0}$ to prevent an unnecessary increase of the marginal probabilities for a part of neighboring points.}
Algorithm~\ref{alg:proposed} shows the pseudo-code of the CMA-ES-SoP.

\subsection{Sample Encoding}
\label{sec:proposed:encoding}
The CMA-ES-SoP transforms the samples $\x^{\langle 1 \rangle}, \cdots, \x^{\langle \lambda \rangle}$ generated from the multivariate Gaussian distribution into the candidate solutions on the search space.
For the elements ${\x}_k^{\langle i \rangle}$ of $\x^{\langle i \rangle}$ corresponding to the $k$-th subspace, the closest points $\textsc{Enc}_k({\x}_k^{\langle i \rangle})$ in ${\S}_k$ are selected as a part of encoded candidate solution, where $\textsc{Enc}_k: \R^{N_k} \to {\S}_k$ is defined as
\begin{align}
    \textsc{Enc}_k({\x}_k) = \argmin_{{\s} \in {\S}_k} \| {\x}_k - {\s} \| \enspace.
    \label{eq:encode}
\end{align}
Then, the encoded candidate solution $\tilde{\x}^{\langle i \rangle} = \textsc{Concat}( \textsc{Enc}_1({\x}_1^{\langle i \rangle}), \cdots, \textsc{Enc}_K({\x}_K^{\langle i \rangle} ))$ is evaluated on the objective function.
The closest point is determined by the Voronoi region containing the corresponding elements of sample ${\x}_k^{\langle i \rangle}$ on the subspace. 

\begin{algorithm*}[!t] 
\centering
\caption{CMA-ES on Sets of Points}
\begin{algorithmic}[1] \label{alg:proposed}
\REQUIRE The objective function $f:\R^N \to \R$ and subspaces $\S_1, \cdots, \S_K$.
\REQUIRE Initial distribution parameters $\mv^{(0)}, \cov^{(0)}, \sigma^{(0)}$.
\REQUIRE Hyperparameters $\alpha_\mathrm{target} = 1 / (\lambda N)$ and $\beta = 1 + 1/N$.
\STATE Initialize the margin value as $\alpha_k^{(0)} = \alpha_\mathrm{target}$.
\WHILE{termination condition is not met}
\FOR{$i = 1, \cdots, \lambda$}
\STATE Generate $\x^{\langle i \rangle}$ from the multivariate Gaussian distribution using (\ref{eq:cma:generation}).
\STATE Encode $\x^{\langle i \rangle}$ to $\tilde{\x}^{\langle i \rangle}$ by concatenating the nearest points as \eqref{eq:encode}.
\STATE Evaluate $\tilde{\x}^{\langle i \rangle}$ on the objective function $f$.
\ENDFOR
\STATE Update $\mv^{(t)}, \cov^{(t)}, \sigma^{(t)}$ using the update rules (\ref{eq:cma:update-mean}), (\ref{eq:cma:update-stepsize}), and (\ref{eq:cma:update-cov}) with samples $\x^{\langle 1 \rangle}, \cdots, \x^{\langle \lambda \rangle}$ before encoding.
\FOR{$k = 1, \cdots, K$}
\STATE Compute the neighboring points $\S_k^\mathrm{neighbor}$ to the mean vector in $k$-th subspace.
\FOR{$\s_{k,b}^\mathrm{neighbor} \in \S_k^\mathrm{neighbor}$}
\STATE Compute the marginal probability $p^{(t+1)}_{k,b}$ in \eqref{eq:margin:p}.
\IF{$p^{(t+1)}_{k,b} < \alpha_k^{(t)}$}
\STATE Correct $\cov^{(t+1)}$ by margin correction with $\alpha_k^{(t)}$ as \eqref{eq:correct:cov}.
\ENDIF
\ENDFOR
\STATE Adjust the margin $\alpha_k^{(t)}$ using the probabilities $p^{(t+1)}_{k,1}, \cdots, p^{(t+1)}_{k,B^{(t+1)}}$ as \eqref{eq:margin:adjust}.
\ENDFOR
\STATE $t \leftarrow t + 1$
\ENDWHILE
\end{algorithmic} 
\end{algorithm*}
\begin{figure}[t]
    \centering
    \begin{minipage}[b]{0.3\linewidth}
        \centering
        \includegraphics[scale=0.3]{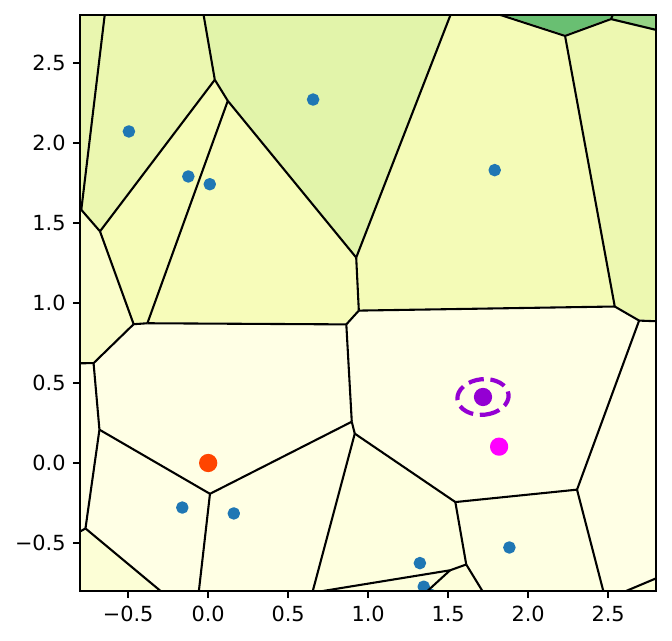}
        \subcaption{Before correction}
    \end{minipage}
    \begin{minipage}[b]{0.3\linewidth}
        \centering
        \includegraphics[scale=0.3]{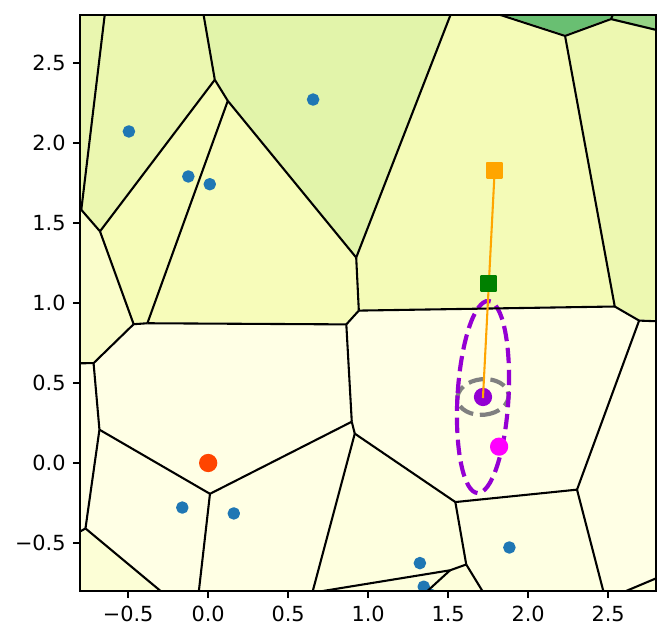}
        \subcaption{A step of correction} 
        \label{fig:margin:B}
    \end{minipage}
    \begin{minipage}[b]{0.3\linewidth}
        \centering
        \includegraphics[scale=0.3]{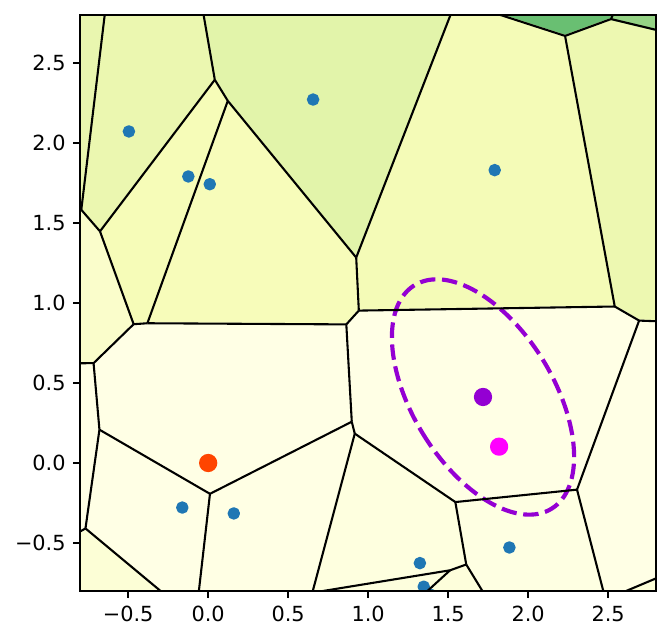}
        \subcaption{After correction}
    \end{minipage}
    \caption{
        Illustration of margin correction in CMA-ES-SoP. The orange and green square points in the center figure represent the neighboring point and mid-point, respectively.
        The gray and violet ellipses in the center figure correspond to the covariance matrices of the multivariate Gaussian distribution before and after a single step of margin correction, respectively.
    }
    \label{fig:margin}
\end{figure}
%
%
%

\subsection{Margin Correction}
The margin correction aims to maintain the generation probability of neighboring points above the margin value to prevent early convergence.
Because the exact computation of generation probability over a non-linearly constrained region is intractable, we develop our margin correction with an alternative tail probability.
Figure~\ref{fig:margin} shows an example of margin correction.

We consider the closest point $\s_k^\mathrm{close} = \textsc{Enc}_k( \mv_k^{(t+1)} )$ to the mean vector on the $k$-th subspace.
We then compute the neighboring points around $\s_k^\mathrm{close}$ on the Voronoi diagram as
\begin{align}
    \S_k^\mathrm{neighbor} = \{ \s_{k, 1}^\mathrm{neighbor}, \cdots, \s_{k, B_k^{(t+1)}}^\mathrm{neighbor} \} \subseteq \S_k \enspace,
\end{align}
\new{where $B_k^{(t+1)}$ is the number of neighboring points.}
With a neighboring point $\s_{k, b}^\mathrm{neighbor}$, the mid-point $\s^\mathrm{mid}_{k, b}$ between the mean vector $\mv_k^{(t+1)}$ and $\s_{k, b}^\mathrm{neighbor}$ is computed as
\begin{align}
     \s^\mathrm{mid}_{k, b} = \frac{\mv_k^{(t+1)} + \s_{k, b}^\mathrm{neighbor}}{2} \enspace.
\end{align}

Then, we consider the marginal distribution along the direction $\s^\mathrm{mid}_{k, b} - \mv_k^{(t+1)}$ \new{to compute the alternative tail probability to exact generation probability}.
We aim to maintain the generation probability $p_{k,b}^{(t+1)}$ beyond the mid-point on the marginal distribution above the margin value $\alpha_k^{(t)}$.
The generation probability is computed as
\begin{align}
    p_{k,b}^{(t+1)} = \Phi_{\mathrm{cdf}} \left( - d_{k,b}^{(t+1)} \right) \enspace,
    \label{eq:margin:p}
\end{align}
where $d_{k,b}^{(t+1)}$ is the Mahalanobis distance between $\s^\mathrm{mid}_{k, b}$ and $\mv_k^{(t+1)}$ on the $k$-th subspace as
\begin{align}
    d_{k,b}^{(t+1)} = \left\| \sqrt{\left( \cov^{(t+1)} \right)^{-1}} \cdot \xi_{k,b}^{(t+1)} \right\| \enspace
\end{align}
with a vector $\xi_{k,b}^{(t+1)}$ defined by two zero vectors $\mathbf{0}_{k}^{\mathrm{ant}}$ and $\mathbf{0}_{k}^{\mathrm{post}}$ with the lengths of $\sum_{j=1}^{k-1} N_j$ and $\sum_{j=k\new{+1}}^{K} N_j$ as
\begin{align}
    \xi_{k,b}^{(t+1)} = \textsc{Concat}\left( \mathbf{0}_{k}^{\mathrm{ant}}, \left( \frac{\s^\mathrm{mid}_{k, b} - \mv_k^{(t+1)}}{\sigma^{(t+1)}} \right), \mathbf{0}_{k}^{\mathrm{post}} \right) \enspace.
\end{align}
When $p_{k,b}$ is smaller than $\alpha_k^{(t)}$, the covariance matrix is modified to
\begin{align}
    \cov^{(t+1)} \leftarrow \cov^{(t+1)} + \frac{ (d_{k,b}^{(t+1)})^2 - (\gamma_\alpha^{(t)})^2 }{ (d_{k,b}^{(t+1)})^2 (\gamma_\alpha^{(t)})^2 } \cdot \xi_{k,b}^{(t+1)} \left( \xi_{k,b}^{(t+1)} \right)^\T \enspace,
    \label{eq:correct:cov}
\end{align}
where $\gamma_\alpha^{(t)} = \Phi_\mathrm{ppf}(1 - \alpha_k^{(t)})$.
This modification maintains the Mahalanobis distance $d_{k,b}^{(t+1)}$ at most $\gamma_\alpha^{(t)}$, which ensures $p_{k,b}^{(t+1)} \geq \alpha_k^{(t)}$ (see appendix~\ref{apdx:sec:ensuring}).
We note that $\xi_{k,b}^{(t+1)}$ can contain a non-zero value on the elements corresponding to $k$-th subspace, and the margin correction does not change the variance on other subspaces.
The CMA-ES-SoP shuffles the neighboring points before each margin correction and applies the correction for each neighboring point in turn.

\subsection{Margin Adaptation}
In the margin correction explained in the previous subsection, a single step in the correction of the covariance matrix in Eq.~\eqref{eq:correct:cov} may increase the marginal probabilities for other neighboring points more than necessary. \del{Therefore}{}\new{To prevent the performance deterioration due to such unnecessary increase of the marginal probabilities}, we adjust the margin value $\alpha_k^{(t)}$ so that the average of probabilities $p_{k,1}^{(t+1)}, \cdots, p_{k,B_k^{(t+1)}}^{(t+1)}$ on marginal distribution is maintained close to the target margin value $\alpha_\mathrm{target}$.
We realize this adjustment by employing the update rule for $\alpha_k^{(t)}$ given by
\begin{align}
    \alpha_k^{(t+1)} = \begin{cases}
    \alpha_k^{(t)} / \beta & \text{if} \quad \alpha_\mathrm{target} \leq \frac{1}{B_k^{(t+1)}} \sum_{b=1}^{B_k^{(t+1)}} p_{k,b}^{(t+1)} \\
    \beta \cdot \alpha_k^{(t)} & \text{otherwise}
    \end{cases}
    \enspace.
    \label{eq:margin:adjust}
\end{align}
We set the target margin value and the increasing and decreasing factor as $\alpha_\mathrm{target} = 1 / (N \lambda)$ and $\beta = 1 + 1/N$, respectively.
We note the target margin value follows the reference~\cite{cmaeswm}.

\begin{figure*}[t]
\centering
\includegraphics[width=0.9\textwidth]{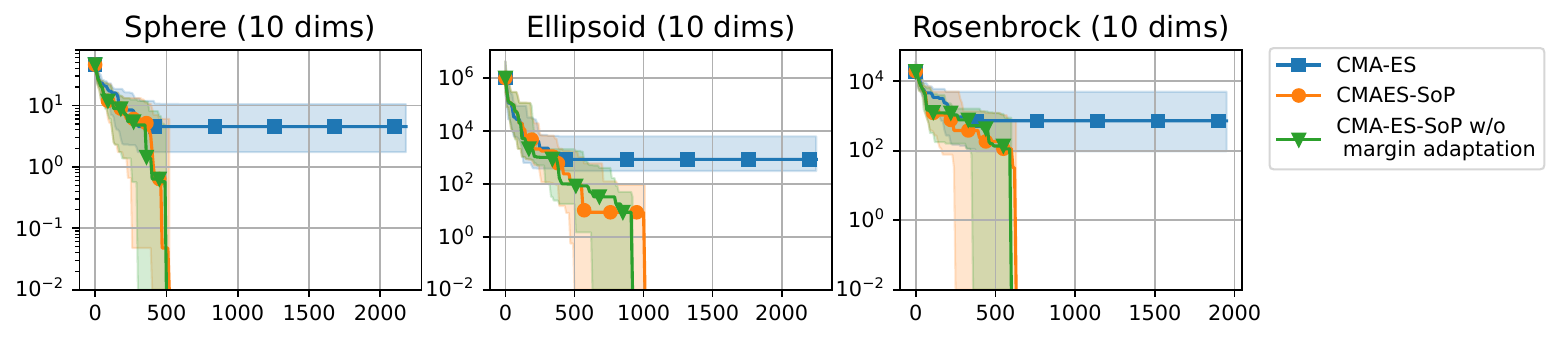}
\includegraphics[width=0.9\textwidth]{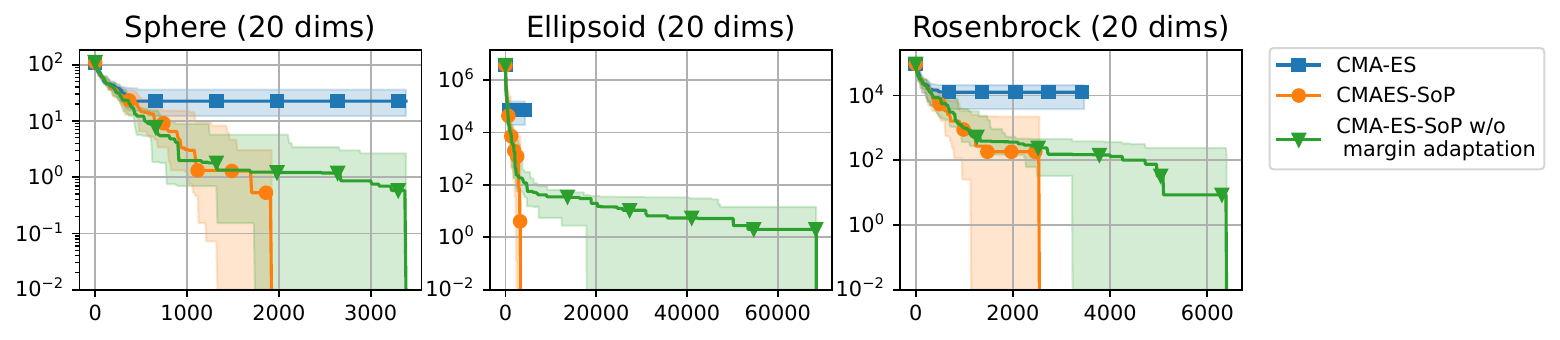}
\includegraphics[width=0.9\textwidth]{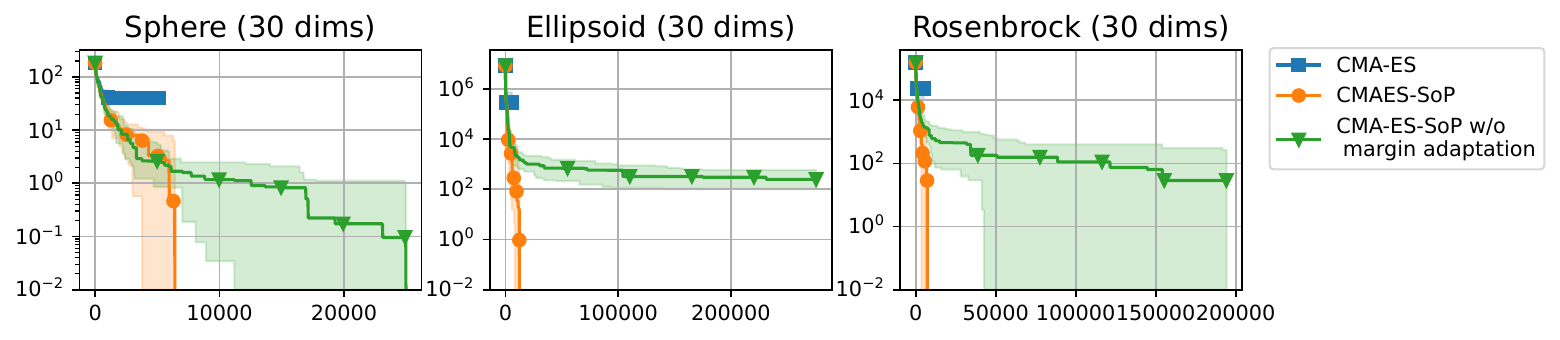}
\caption{
Transitions of the best evaluation values on the discrete optimization problems with $(N_k, L_k) = (2,10)$.
We plot the median and interquartile ranges over 25 independent trials.
}
\label{fig:exp1}
\end{figure*}
\begin{table}[t]
  \caption{Success rate (SR) and SP1 in discrete optimization on sets of points.}
  \label{table:discrete}
  \centering
  \vspace{5pt}
  \begin{tabular}{c|c|c||c|c||c|c||c|c|}
    \multicolumn{2}{c|}{\multirow{2}{*}{Problem Setting}} & \multirow{2}{*}{Method}& \multicolumn{2}{c||}{Sphere} & \multicolumn{2}{c||}{Ellipsoid} & \multicolumn{2}{c|}{Rosenbrock} \\
    \multicolumn{2}{c|}{} & & \hspace{4pt} SR \hspace{4pt} & SP1 & \hspace{4pt} SR \hspace{4pt} & SP1 & \hspace{5pt} SR \hspace{5pt} & SP1 \\
    \hline
    \hline
    \multirow{6}{4em}{\centering $N_k=2$ \par $L_k=10$} & \multirow{2}{*}{\hspace{1pt} $N = 10$ \hspace{1pt}} & CMA-ES & 0.20 & {\bf 1410.0} & 0.00 & -- & 0.24& {\bf 1069.4} \\
    & & CMA-ES-SoP & {\bf 1.00} & 1611.2 & {\bf 0.96} &  {\bf 1406.6} & {\bf 0.96} & 1282.1 \\
    \cline{2-9}
    & \multirow{2}{*}{$N = 20$} & CMA-ES & 0.00 & -- & 0.00 & -- & 0.00 & -- \\
    & & CMA-ES-SoP & {\bf 1.00} & {\bf 3811.6} & {\bf 1.00} & {\bf 5002.5} & {\bf 1.00} & {\bf 6043.6} \\
    \cline{2-9}
    & \multirow{2}{*}{$N = 30$} & CMA-ES & 0.00 & -- & 0.00 & -- & 0.00 & -- \\
    & & CMA-ES-SoP & {\bf 1.00} & {\bf 9456.1} & {\bf 1.00} & {\bf 12291.4} & {\bf 0.96} & {\bf 12534.9} \\
    \hline
    \hline
    \multirow{6}{4em}{\centering $N_k=5$ \par $L_k=40$} & \multirow{2}{*}{$N = 10$} & CMA-ES & 0.32 & 277.3 & 0.12 & 805.5 & 0.44 & {\bf 78.5} \\
    & & CMA-ES-SoP & {\bf 1.00} & {\bf 213.2}  & {\bf 1.00} & {\bf 541.6} & {\bf 1.00} & 134.8 \\
    \cline{2-9}
    & \multirow{2}{*}{$N = 20$} & CMA-ES & 0.00 & -- & 0.00 & -- & 0.04 & 4800.0 \\
    & & CMA-ES-SoP & {\bf 1.00} & {\bf 765.6} & {\bf 1.00} & {\bf 4431.3} & {\bf 0.96} & {\bf 1679.6} \\
    \cline{2-9}
    & \multirow{2}{*}{$N = 30$} & CMA-ES & 0.00 & -- & 0.00 & -- & 0.00 & -- \\
    & & CMA-ES-SoP & {\bf 1.00} & {\bf 2107.28} & {\bf 1.00} & {\bf 7458.6} & {\bf 1.00} & {\bf 2667.2} \\
    \hline
  \end{tabular}
\end{table}
\begin{figure*}[t]
\centering
\includegraphics[width=0.9\textwidth]{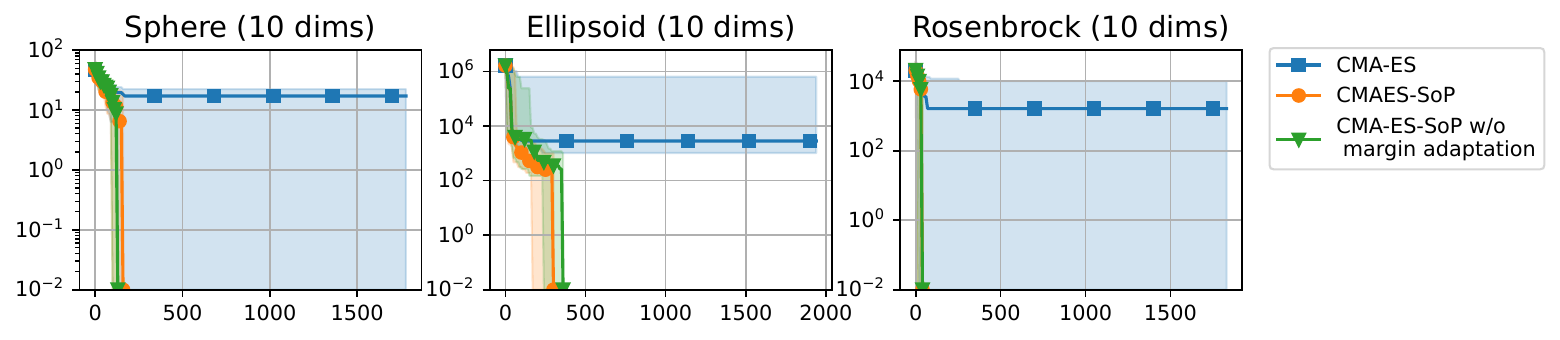}
\includegraphics[width=0.9\textwidth]{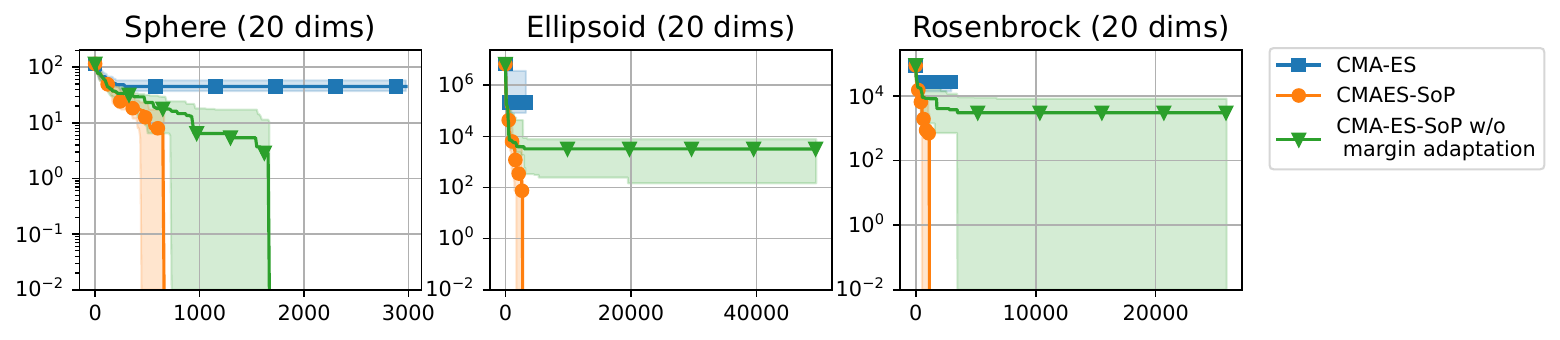}
\includegraphics[width=0.9\textwidth]{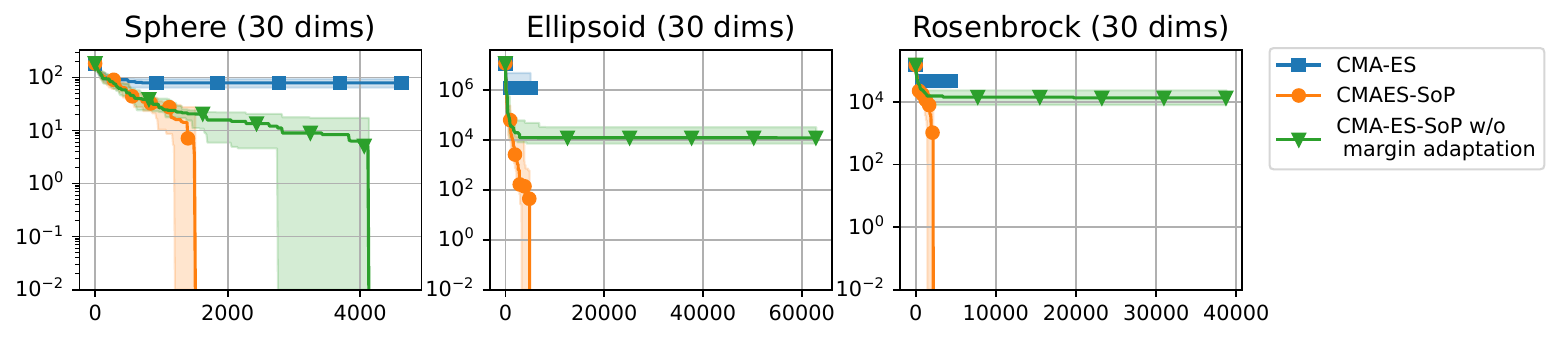}
\caption{
Transitions of the best evaluation values on the discrete optimization problems with $(N_k, L_k) = (5,40)$.
We plot the median and interquartile ranges over 25 independent trials.
}
\label{fig:exp2}
\end{figure*}
\begin{figure*}[t]
\centering
\includegraphics[width=0.8\textwidth]{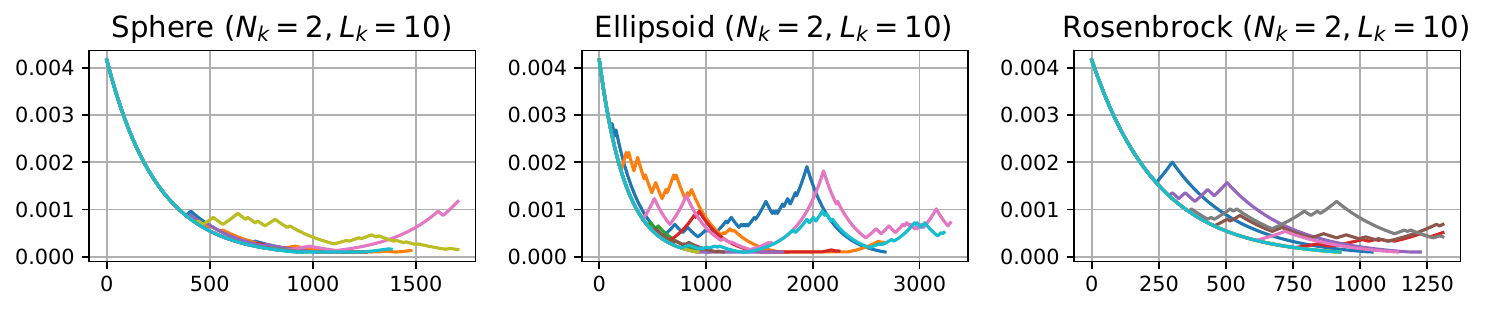}
\includegraphics[width=0.8\textwidth]{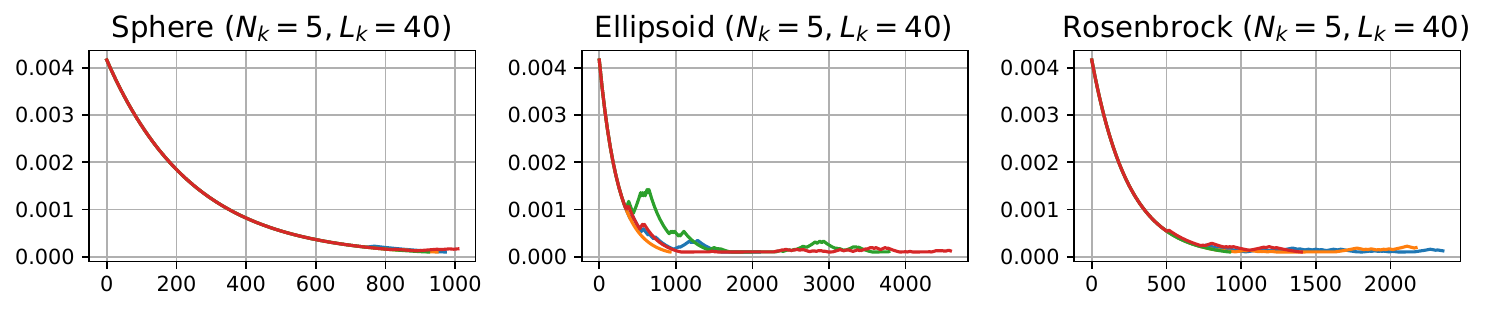}
\caption{
Transitions of margins for each subspace. They were observed in a typical trial of the CMA-ES-SoP on discrete 20-dimensional optimization problems.
}
\label{fig:margin:disc}
\end{figure*}
%
%
%

\section{Experiment}
\label{sec:exp}
We evaluated the optimization performance of the CMA-ES-SoP on the discrete optimization problems on the sets of points in Section~\ref{sec:exp:disc} and the mixed-variable optimization problems in Section~\ref{sec:exp:mixed}.
The code of the CMA-ES-SoP will be made available at \textcolor{blue}{\url{https://github.com/CyberAgentAILab/cmaes}}~\cite{nomura2024cmaes}.

\subsection{Experimental Setting}
We prepared four benchmark functions as follows:
\begin{itemize}
    \item Sphere: $f(\x) = \sum_{i=1}^{N} x_i^2$
    \item Ellipsoid: $f(\x) = \sum_{i=1}^{N} \left( 1000^{\frac{i-1}{{N}-1}} x_i \right)^2$
    \item Reversed Ellipsoid: $f(\x) = \sum_{i=1}^{N} \left( 1000^{\frac{{N}-i}{{N}-1}} x_i \right)^2$
    \item Rosenbrock: $f(\x) = \sum_{i=1}^{{N}-1} \left( 100 (x_{i+1}  - x_i^2)^2 + (x_i - 1)^2 \right)$
\end{itemize}
\rev{Sphere, Ellipsoid, and Rosenbrock are well-known benchmark functions. We added Reversed Ellipsoid for comparison with Ellipsoid, which deepens our discussions in mixed-variable optimization. 
}

The sets of points $\S_k$ was given by $L_k - 1$ samples generated from the uniform distribution on $[-5, 5]^{N_k}$.
We then added the optimal solution of benchmark functions to $\S_k$.
We performed two settings for $N_k$ and $L_k$ as $(N_k, L_k) \in \{ (2, 10), (5, 40) \}$.
We varied the total number of dimensions as $N = 10, 20, 30$.

We ran the naive CMA-ES explained in Section~\ref{sec:cmaes} \new{(with the sample encoding in Section~\ref{sec:proposed:encoding})} as a comparative method. 
We also compared the CMA-ES-SoP without the margin adaptation, \del{which fixed $\alpha^{(t)}_k = \alpha_\mathrm{target}$}{}\new{in which $\alpha^{(t)}_k = \alpha_\mathrm{target}$ was fixed} for all iterations.
For both CMA-ES and CMA-ES-SoP, the initial mean vector $\mv^{(0)}$ was given uniformly at random on $[1,5]^N$.
The initial covariance and step-size were given by $\cov^{(0)} = \I$ and $\sigma^{(0)} = 2$, respectively.
We terminated the optimization when one of the following four conditions was met: 1) the successful condition was satisfied, 2) the number of evaluations reached $N \times 10^4$, 3) the minimum eigenvalue of $(\sigma^{(t)})^2 \cov^{(t)}$ was updated less than $10^{-30}$, or 4) a numerical error occurred.


\subsection{Experimental Result in Discrete Optimization on Sets of Points}
\label{sec:exp:disc}
First, we show the results with the discrete optimization problems on sets of points.
In this setting, the search space consisted of $K = N / N_k$ sets of points.
We regarded a trial as successful when the optimal solution was found.
We ran 25 independent trials on Sphere, Ellipsoid, and Rosenbrock for each setting.

Figures~\ref{fig:exp1} and~\ref{fig:exp2} show the transitions of the best evaluation values with $(N_k, L_k) = (2,10)$ and $(N_k, L_k) = (5,40)$, respectively.
We observed that the CMA-ES-SoP successfully optimized all benchmark problems while the CMA-ES \rev{often} stagnated on all benchmark functions.
Note that the optimization of the CMA-ES is usually terminated due to a too small eigenvalue of the covariance matrix on the search distribution.
We consider that the margin correction in CMA-ES-SoP prevented such premature convergence.
We also observed that margin adaptation improved the optimization performance for high-dimensional problems.
Figure~\ref{fig:margin:disc} shows the transition of the margins in a typical trial.
We can see that the dynamics of the margin change depending on the function, subspace dimension, and number of data, highlighting the importance of margin adaptation.

Table~\ref{table:discrete} shows the success rates and SP1 values computed with 25 trials.
SP1 is the average number of evaluations over successful trials divided by the success rate~\cite{sp1}.
We can see that the success rates of the CMA-ES-SoP were significantly better than those of the CMA-ES in all problem settings.
Focusing on SP1, however, the CMA-ES was superior to the CMA-ES-SoP on some low-dimensional problems. \revdel{We consider that the CMA-ES quickly converged to the optimum solution when the initial mean vector was given close to the optimal solution.}\rev{We note that, although the CMA-ES failed to optimize in most of trials, it sometimes quickly converged to the optimum solution.}

\begin{figure*}[t]
\centering
\includegraphics[width=0.9\textwidth]{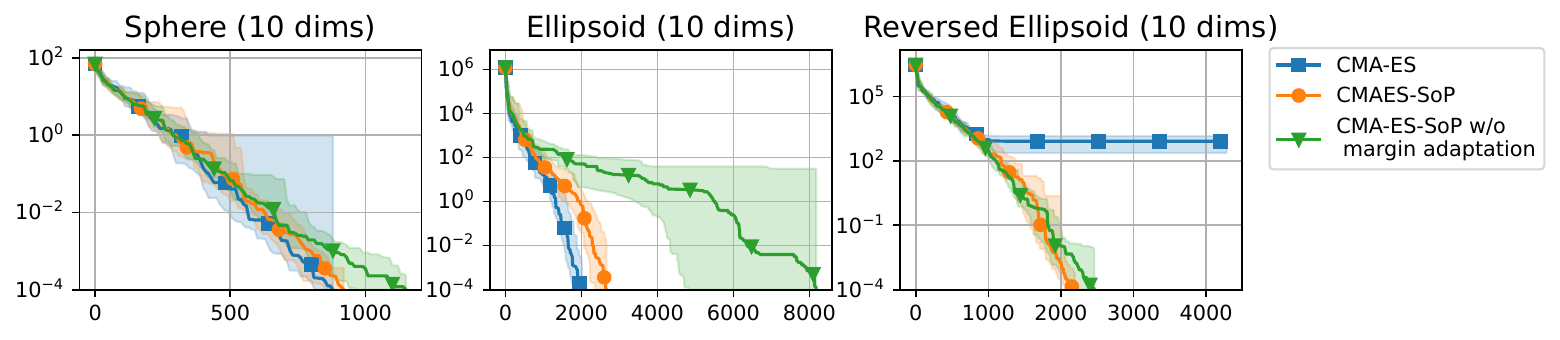}
\includegraphics[width=0.9\textwidth]{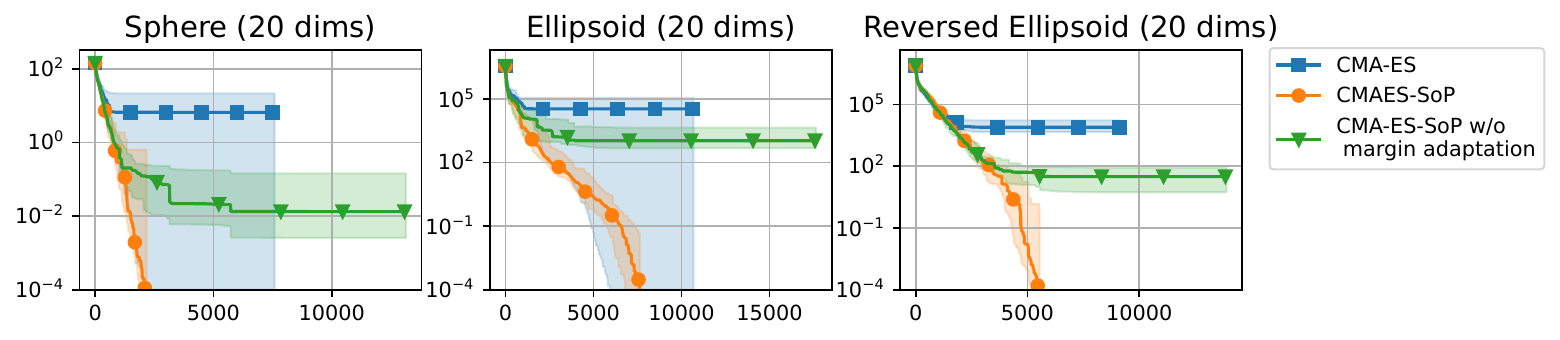}
\includegraphics[width=0.9\textwidth]{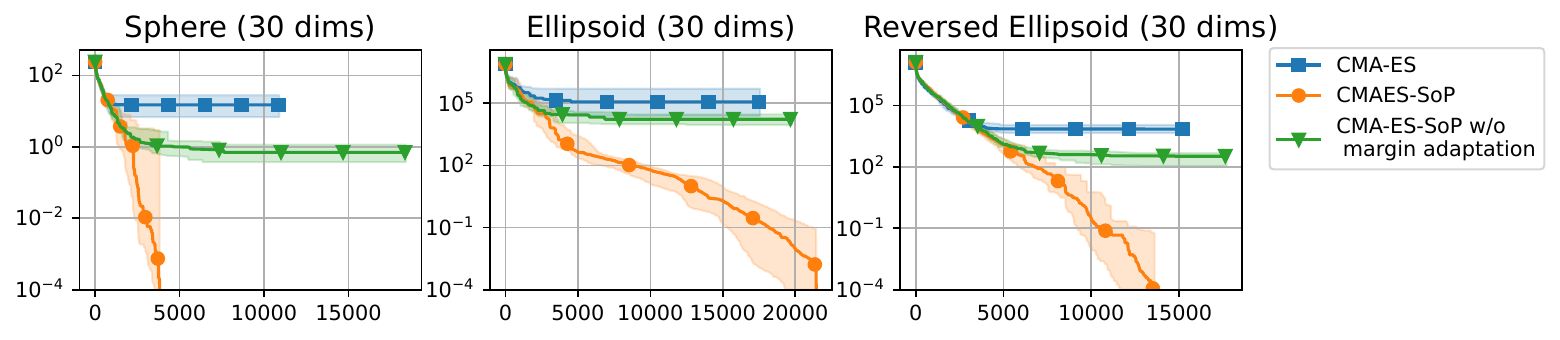}
\caption{
Transitions of the best evaluation values on the mixed-variable optimization problems with $(N_k, L_k) = (2,10)$.
We plot the median and interquartile ranges over 25 independent trials.
}
\label{fig:exp4}
\end{figure*}
\begin{figure*}[t]
\centering
\includegraphics[width=0.9\textwidth]{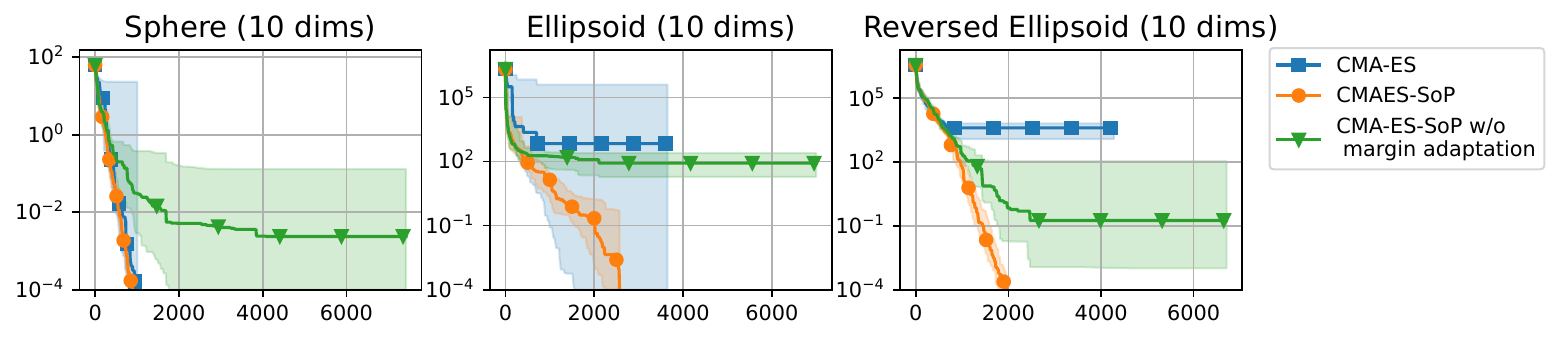}
\includegraphics[width=0.9\textwidth]{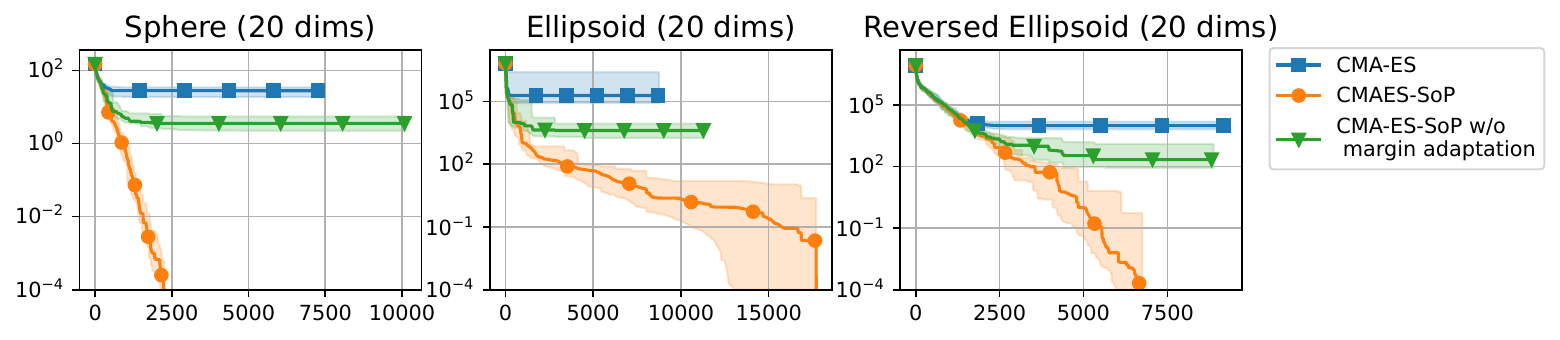}
\includegraphics[width=0.9\textwidth]{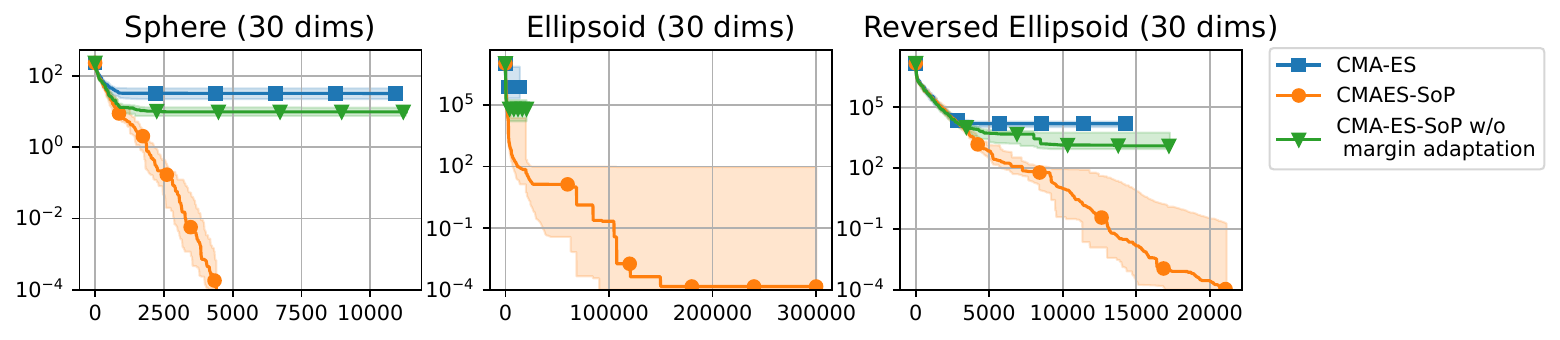}
\caption{
Transitions of the best evaluation values on the mixed-variable optimization problems with $(N_k, L_k) = (5,40)$.
We plot the median and interquartile ranges over 25 independent trials.
}
\label{fig:exp5}
\end{figure*}
\begin{table}[t]
  \caption{Success rate (SR) and SP1 in mixed-variable optimization.}
  \label{table:mixed}
  \centering
  \vspace{5pt}
  \begin{tabular}{c|c|c||c|c||c|c||c|c|}
    \multicolumn{2}{c|}{} && \multicolumn{2}{c||}{\multirow{2}{*}{Sphere}} & \multicolumn{2}{c||}{\multirow{2}{*}{Ellipsoid}} & \multicolumn{2}{|c|}{Reversed} \\
    \multicolumn{2}{c|}{Problem Setting} & Method & \multicolumn{2}{c||}{} & \multicolumn{2}{c||}{} & \multicolumn{2}{c|}{Ellipsoid}\\
    \multicolumn{2}{c|}{} && \hspace{5pt}SR\hspace{5pt} & SP1 & \hspace{5pt}SR\hspace{5pt} & SP1 & \hspace{5pt}SR\hspace{5pt} & SP1 \\
    \hline
    \hline
    \multirow{6}{4em}{\centering $N_k=2$ \par $L_k=10$} & \multirow{2}{*}{\hspace{1pt}$N = 10$\hspace{1pt}} & CMA-ES & 0.72 & 2120.3 & 0.84 & {\bf 3075.9} & 0.00 & -- \\
    & & CMA-ES-SoP & {\bf 1.00} & {\bf 1567.2} & {\bf 1.00} & 3652.8 & {\bf 1.00} & {\bf 3605.6}  \\
    \cline{2-9}
    & \multirow{2}{*}{$N = 20$} & CMA-ES & 0.28 & 11448.9 & 0.28 & 22646.9 & 0.00 & -- \\
    & & CMA-ES-SoP & {\bf 1.00} & {\bf 3632.6} & {\bf 1.00} & {\bf 10402.5} & {\bf 1.00} & {\bf 14764.8} \\
    \cline{2-9}
    & \multirow{2}{*}{$N = 30$} & CMA-ES & 0.04 & 120750.0 & 0.12 & 101188.8 & 0.00 & -- \\
    & & CMA-ES-SoP & {\bf 1.00} & {\bf 6444.4} & {\bf 1.00} & {\bf 25319.2} & {\bf 1.00} & {\bf 26569.7} \\
    \hline
    \hline
    \multirow{6}{4em}{\centering $N_k=5$ \par $L_k=40$} & \multirow{2}{*}{$N = 10$} & CMA-ES & 0.52 & 2871.3 & 0.40 & {\bf 5455.0} & 0.08 & 29687.5 \\
    & & CMA-ES-SoP & {\bf 1.00} & {\bf 1594.0} & {\bf 0.92} & 12787.3 & {\bf 0.92} & {\bf 7545.3}  \\
    \cline{2-9}
    & \multirow{2}{*}{$N = 20$} & CMA-ES & 0.04 & 19837.5 & 0.04 & {\bf 44300.0} & 0.00 & -- \\
    & & CMA-ES-SoP & {\bf 0.96} & {\bf 3835.4} & {\bf 0.76} & 78968.1 & {\bf 0.84} & {\bf 57559.8} \\
    \cline{2-9}
    & \multirow{2}{*}{$N = 30$} & CMA-ES & 0.04 & 119700.0 & 0.12 & {\bf 238000.0} & 0.00 & -- \\
    & & CMA-ES-SoP & {\bf 1.00} & {\bf 6890.8} & {\bf 0.48} & 355264.5 & {\bf 0.64} & {\bf 185078.9} \\
    \hline
  \end{tabular}
\end{table}
%
%
%

\subsection{Experimental Result in Mixed-Variable Optimization}
\label{sec:exp:mixed}
Next, we show the results with the mixed-variable optimization problems.
In this setting, we set the number of sets of points as \new{$N_\mathrm{set} = \lfloor N / N_k / 2 \rfloor$}.
In this experimental setting, the search space of the first $N_\mathrm{set} N_k$ design variables was given by the sets of points, while the remaining \del{$N - N_\mathrm{set} N_k$}{}\new{$N_K = N - N_\mathrm{set} N_k$} design variables were treated as continuous variables.
We regarded a trial as successful when the best evaluation value reached $10^{-4}$.
We ran 25 independent trials on Sphere, Ellipsoid, and Reversed Ellipsoid for each setting.

Figures~\ref{fig:exp4} and~\ref{fig:exp5} show the transitions of the best evaluation values with $(N_k, L_k) = (2,10)$ and $(N_k, L_k) = (5,40)$, respectively.
On 20- and 30-dimensional problems, the CMA-ES-SoP outperformed the CMA-ES.
In contrast, the CMA-ES was sometimes competitive or superior to the CMA-ES-SoP on the 10-dimensional problems.
We note the number of possible points for discrete variables was significantly low \revdel{in}\rev{on} the 10-dimensional problems, which was \revdel{$10^2$}\rev{100} when $(N_k, L_k) = (2, 10)$ and $40$ when $(N_k, L_k) = (5, 40)$.
This is why the CMA-ES could find the optimum points on the sets of points and showed competitive and superior performance compared to the CMA-ES-SoP.
As well as the result in the previous subsection, the margin adaptation improved the optimization performance on high-dimensional problems.
Because we did not tune the increasing and decreasing factor $\beta$, tuning it may improve the optimization performance.

Table~\ref{table:mixed} shows the success rates and SP1 values computed with 25 trials.
In addition to the discrete optimization in the previous subsection, we can see that the success rates of CMA-ES-SoP were significantly higher compared to CMA-ES in all problem settings.
In addition, the CMA-ES-SoP achieved smaller SP1 values than that of CMA-ES under all settings on Sphere and Reversed Ellipsoid.
On Ellipsoid, however, the CMA-ES was sometimes better than the CMA-ES-SoP with respect to SP1.
We consider the following reason: when optimizing the Ellipsoid by the CMA-ES, the variance corresponding to the continuous variables converged faster than the variance corresponding to the discrete variables.
Therefore, unlike on the other functions, the premature convergence of the CMA-ES in discrete subspace was relatively prevented.

\section{Conclusion}

We have proposed an extension of the CMA-ES for discrete and mixed-variable optimization problems on sets of points.
The proposed CMA-ES-SoP contains three additional steps: sample encoding, margin correction, and margin adaptation.
In the sample encoding, the samples generated from multivariate Gaussian distribution are mapped to the closest points in the subspaces.
In the margin correction, the updated covariance matrix $\cov^{(t+1)}$ is modified to maintain the marginal generation probability $p_{k,b}^{(t+1)}$ above the margin $\alpha^{(t)}_k$.
Finally, in the margin adaptation, the margin is adjusted so that the average of marginal probabilities is maintained close to the target value $\alpha_\mathrm{target}$.
The numerical simulation showed the efficiency of CMA-ES-SoP in discrete and mixed-variable optimizations on sets of points.

\revdel{Differently from the CMA-ES with margin and DX-NES-ICI, we only modified the covariance matrix $\cov^{(t+1)}$ and fixed the mean vector in the margin correction.
The development of efficient modification of the mean vector is one of our future works.}
\rev{In our experiment, we used the benchmark functions extended from the benchmarks for continuous optimization. As the benchmark functions for optimization on sets of points are not well-structured, developing suitable benchmark functions is one of our future works.}
Moreover, because the step-size adaptation in the CMA-ES assumes the optimization in continuous space, we will develop an efficient step-size adaptation for optimization on sets of points in the future.



\appendix

\section{Ensuring Margin by Modification of Covariance}
\label{apdx:sec:ensuring}
In the following, we omit the iterators for short, e.g. we denote $\xi_{k,b}^{(t+1)}$ as $\xi_{k,b}$.
According to the Sherman–Morrison formula, the inverse matrix of modified covariance matrix in Eq.~\eqref{eq:correct:cov} is given by
\begin{align}
    \cov^{-1} = \bar{\cov}^{-1} - \frac{ \zeta \cdot \bar{\cov}^{-1} \xi_{k,b} \xi_{k,b}^\T \bar{\cov}^{-1} }{1 + \zeta \cdot \xi_{k,b}^\T  \bar{\cov}^{-1} \xi_{k,b} }
    \enspace,
\end{align}
where $\bar{\cov}$ is the covariance matrix before the margin correction, and $\zeta = (d_{k,b}^2 - \gamma_\alpha^2) / d_{k,b}^2 / \gamma_\alpha^2 $.
Considering the relation $d_{k,b}^2 = \xi_{k,b}^\T \bar{\cov}^{-1} \xi_{k,b}$, the squared Mahalanobis distance after the margin correction is given by
\begin{align}
    \xi_{k,b}^\T \cov^{-1} \xi_{k,b} &= \xi_{k,b}^\T \bar{\cov}^{-1} \xi_{k,b} - \frac{ \zeta \cdot \left( \xi_{k,b}^\T \bar{\cov}^{-1} \xi_{k,b} \right)^2  }{1 + \zeta \cdot \xi_{k,b}^\T  \bar{\cov}^{-1} \xi_{k,b} } = d_{k,b}^2 - \frac{\zeta \cdot d_{k,b}^4}{1 + \zeta \cdot d_{k,b}^2} \enspace.
\end{align}
Then, substituting $\zeta = (d_{k,b}^2 - \gamma_\alpha^2) / d_{k,b}^2 / \gamma_\alpha^2 $ shows $\xi_{k,b}^\T \cov^{-1} \xi_{k,b} = \gamma_\alpha^2 $.
Finally, remaining that $\Phi_\mathrm{ppf}$ is the inverse function of $\Phi_\mathrm{cdf}$ and $\gamma_\alpha = \Phi_\mathrm{ppf}(1 - \alpha_k) = - \Phi_\mathrm{ppf}(\alpha_k)$, we have $p_{k,b} = \alpha_k$ after the margin correction.

\bibliographystyle{splncs04}
\bibliography{reference}

\end{document}